\input amstex

\documentstyle{amsppt}
\hsize = 5.4 truein
\vsize = 8.7 truein
\NoBlackBoxes
\nologo
\TagsAsMath
\baselineskip = 14pt
\topmatter
\title Coding and Reshaping \\
When There Are No Sharps
\endtitle
\author Saharon Shelah$^{1,2,5}$ and Lee J. Stanley$^{3,4,5}$
\endauthor
\thanks {\roster\item"{1.}" Research partially supported by NSF,
the Basic Research Fund of the Israel Academy of Science, and MSRI.\item"{2.}"
Paper number 294.\item"{3.}"
Research partially supported by NSF grant DMS 8806536 and MSRI.
\item"{4.}" Preparation of the final
version of this paper partially supported by a grant from the Reidler
Foundation.\item"{5.}" We are grateful to the administration and staff
of MSRI for their hospitality during portions of 1989-90.\endroster}\endthanks
\address{Hebrew University \& Rutgers University}\endaddress
\address{Lehigh University}\endaddress
\abstract{Assuming $0^\sharp$ does not exist, $\kappa$ is an uncountable
cardinal and for all cardinals $\lambda$ with $\kappa \leq \lambda <
\kappa^{+\omega},\ 2^\lambda = \lambda^+$,
we present a \lq\lq mini-coding" between $\kappa$ and $\kappa^{+\omega}$.
This allows us to prove that any subset of $\kappa^{+\omega}$ can be coded
into a subset, $W$ of $\kappa^+$ which, further, \lq\lq reshapes" the interval
$[\kappa,\ \kappa^+)$, i.e., for all $\kappa < \delta < \kappa^+,
\ \kappa = (card\ \delta)^{L[W \cap \delta]}$.  We sketch two applications of
this result, assuming $0^\sharp$ does not exist.  First, we point out that
this shows that any set can be coded by a real, via a \bf set \rm forcing.
The second application involves a notion of abstract condensation, due to
Woodin.  Our methods can be used to show that for any cardinal $\mu$,
condensation for $\mu$ holds in a generic extension by a \bf set \rm forcing.}
\endabstract
\endtopmatter
\subheading{\S 0. INTRODUCTION}
\bigskip
\proclaim{Theorem}  Assume that $V \models ZFC + \lq\lq 0^\sharp$ does not
exist", and, in $V,\ \kappa \geq \aleph_2,
\ Z \subseteq \kappa^{+\omega}$ and
for cardinals $\lambda$ with $\kappa \leq \lambda < \kappa^{+\omega},
\ 2^\lambda = \lambda^+$.
\bf THEN \rm there is a cofinality preserving forcing $\bold {S}(\kappa)
= \bold {S}(\kappa,\ Z)$ of cardinality $\kappa^{+(\omega+1)}$
such that if $G$ is $V$-generic for $\bold {S}(\kappa)$, there is
$W \subseteq \kappa^+$ such that $V[G] = V[W],
\ Z \in L[W,\ Z \cap \kappa]$, for all
cardinals $\lambda$ with $\kappa \leq \lambda < \kappa^{+\omega}$, and
for all limit ordinals $\delta$ with $\kappa < \delta < \kappa^+,
\ \kappa = (card\ \delta)^{L[W \cap \delta]}$.
\endproclaim
\bigskip
Our forcing $\bold {S}(\kappa)$ can be thought of as a kind of
Easton product between $\kappa$ and $\kappa^{+\omega}$ of
partial orderings which \it simultaneously \rm
perform the tasks of coding (\S 1.2 of \cite {1}) and reshaping (\S 1.3
of \cite {1}).  Our new idea is to introduce an additional coding
area used for \lq\lq marking" certain ordinals.
This \lq\lq marking" technique is the crucial addition to the arguments
of \S 1 of \cite{1}.  We appeal to the Covering Lemma twice:  in (2.1), and
again in the proof of the Proposition in (2.3).
The referee has informed us that the hypothesis that $0^\sharp$ does not
cannot be eliminated.  Jensen first used this hypothesis in \cite{1}
to facilitate certain arguments, and then realized that {\it his} uses
{\it were} eliminable.  It is not the purpose of this paper to discuss
the nature of Jensen's appeals to the Covering Lemma; the
interested reader may consult pp. 62, 96 and the Introduction to Chapter 8
of \cite{1} for insight into Jensen's uses of the Covering Lemma, and how he
was able to eliminate them.  In \cite{2}, S. Friedman presents a rather
different, more streamlined approach to avoiding such uses of Covering.
It should be clear from
the preceding that Jensen's appeals to the Covering Lemma
are of a rather different character than ours.

To better understand the role of this \lq\lq marking" technique, let
us briefly recall some material from \cite{1}.  Let us first consider
the possibility of coding
$R \subseteq \kappa^+$ into a subset of
$\kappa$, when $\kappa$ is regular.  In order to use
almost disjoint set coding,
we seem to need extra properties of the ground
model, or of the set $R$, since, in order to carry out the {\it decoding}
recursion across $[\kappa,\ \kappa^+)$ we need, e.g., an almost disjoint
sequence $\overset{\rightharpoonup}\to{b} =
(b_\alpha:\alpha\in [\kappa,\ \kappa^+))$ of cofinal subsets of $\kappa$
satisfying:
$$
\split
(\ast):\ \ \ \text{  for all } \theta \in [\kappa,\ \kappa^+),
\ (b_\alpha: \alpha \leq \theta) \in L[R \cap \theta],\\
\text{ and is \lq\lq canonically definable" there}.
\endsplit
$$
Such a $\overset{\rightharpoonup}\to{b}$ is called
{\it decodable}, and it is easy to obtain a  decodable
$\overset{\rightharpoonup}\to{b}$ if $R$ satisfies:
$$
(\ast\ast):\text{  for all  } \theta \in [\kappa,\ \kappa^+),
\ (card\ \theta)^{L[R \cap \theta]}= \kappa.
$$
If ($\ast\ast$) holds, we say that $R$ {\it promptly collapses fake cardinals}.

Of course, typically ($\ast\ast$) fails, and the
\lq\lq reshaping" conditions of \S 1.3 of
\cite{1} are introduced
to obtain ($\ast\ast$) in a generic
extension.  Our $\kappa$ and $R$, from the previous paragraph
are called $\gamma$ and $B$ in \S 1.3 of \cite{1}.
Unfortunately, the distributivity argument for the reshaping
partial ordering given there {\it seems}
to really require not merely that $H_{\gamma^+} = L_{\gamma^+}[B]$,
but that $H_{\gamma^{++}} = L_{\gamma^{++}}[B]$,
where $B \subseteq
\gamma^+$.  This will be the case if $B$ is the result of coding as far as
$\gamma^{+}$, but that is another story,
which leads to Jensen's original approach to
the Coding Theorem.  Our appeals to the Covering Lemma focus on this point:
essentially, to prove a distributivity property of the reshaping conditions.
As already indicated, in Jensen's treatment, the appeals to the Covering Lemma
were designed to overcome different sorts of obstacles and proved to be
eliminable.

Our approach to guaranteeing that the unions of certain increasing chains of
reshaping conditions collapse the suprema of their domains is to have
\lq\lq marked" a cofinal sequence of small order type.
Because of the need to meet certain dense sets in the
course of the construction, it is too much to expect that the ordinals
we \bf intentionally \rm marked are the only marked ordinals.  However,
what we will be able to guarantee is that they are the only members
of a certain club subset which have been marked.  The club will exist in
a small enough inner model, thanks to the Covering Lemma.
This argument is given in (2.3).
We are grateful to the referee for suggesting
the use of \lq\lq fast clubs" in the argument of (2.3).  This allowed us to
streamline a more complicated argument (which also suffered from some
[probably reparable] inaccuracies) in an earlier version of this paper.
We use \lq\lq 1" to mark ordinals.  To guarantee that this does not collide
with requirements imposed by the \lq\lq coding" part of the conditions, we set
aside the limit ordinals as the only potentially marked ordinals and do not use
them for coding.

\bigskip

\noindent
SUMMARY AND ORGANIZATION.
\medskip

We now give a brief overview of the contents of this paper.  In \S 1, we
build to the definition, in (1.5), of the $\bold {S}(\kappa)$, along with
auxiliary forcings, $\bold {S}_k(\kappa)$.  In \S 2, we prove that the
$\bold {S}(\kappa)$ are as required.  The heart of the matter is
(2.3), where we prove the distributivity properties of the
$\bold {S}_k(\kappa)$.  Preliminary observations are given in (2.1) and (2.2).
The former shows that only increasing sequences of certain lengths are
problematical.  The latter is a rather routine observation about how the
coding works.  In the argument of (2.3), we use this in the context of forcing
over $\hat{\Cal N}$, a transitive set model of enough $ZFC$,
introduced in the proof of (2.3), below.
In (2.4) we put together the material of (2.1) - (2.3)
to prove the Theorem.  In (2.5) we make a few remarks and briefly sketch the
applications mentioned in the abstract.

The partial ordering $\bold {S}(T,\ \lambda)$,
introduced in (1.2), below, is the analogue of the reshaping partial ordering
of \S (1.3) of \cite{1}.  It adds a subset of $\lambda^+$, which, together with
$T$, promptly collapses fake cardinals in $(\lambda,\ \lambda^+)$.
The partial
ordering $\bold {P}_{\kappa,\ T,\ g}$, introduced in (1.4), is a version
of the coding partial ordering of (1.2) of \cite{1}, \bf relative to \rm $g$.
We require that $T \subseteq \kappa^+,\ g \in S(T,\ \kappa^+)$.
If $p \in P_{\kappa,\ T,\ g}$, then $p$ will have the form
$(\ell(p),\ r(p));\ \ell(p)$ is the \lq\lq function part" of $p$ and
$r(p)$ is the \lq\lq promise part" of $p$.  We require that
$\ell(p)$ starts to code not only $T$, but also $g$ \bf and \rm
that $\ell(p) \in S(T \cap \kappa,\ \kappa)$.  If $g$ were not merely
a condition but generic for $\bold{S}(T,\ \kappa^+)$, then
$\bold {P}_{\kappa,\ T,\ g}$ would just be the usual forcing for the
almost-disjoint set coding of the \lq\lq join" of $T$ and $g$, with the extra
requirement above, that for conditions, $p,\ \ell(p)$, together with
$T \cap \kappa$, collapses $sup\ dom\ \ell(p)$.

Finally, the $\bold {S}(\kappa)$, introduced in (1.5), is the forcing which
accomplishes the task of coding and reshaping, between $\kappa$ and
$\kappa^+$.  It is defined relative to the choice of
a fixed $Z \subseteq \kappa^{+\omega}$
such that $H_{\kappa^{+n}} =
L_{\kappa^{+n}}[Z \cap \kappa^{+n}]$, for all
$n \leq \omega$.  The elements of $S(\kappa)$, are $\omega\text{-sequences},
\ (p(n):\ n < \omega)$, where for all $n < \omega,
\ p(n) = (\ell(p(n)),\ r(p(n))),\ \ell(p(n)) \in
S(Z \cap \kappa_n,\ \kappa_n)$ and
$p(n) \in  P_{\kappa_n,\ Z \cap \kappa_{n+1},\ \ell(p(n+1))}$.
Thus, letting
$\dot G$  be the canonical name
for the generic of $\bold {S}(\kappa)$,
letting $\dot {G}(n)$ be the canonical name for
$\{\ell(p(n)): p \in \dot G\}$, and letting $\dot {W}(n)$ be the
canonical name for $\bigcup \dot {G}(n)$,
$\bold {S}(\kappa)$ is
a sort of Easton product of the
$P_{\kappa_n,\ Z \cap \kappa_{n+1},\ \dot {W}(n+1)}$.
\bigskip

\noindent
NOTATION AND TERMINOLOGY.
\medskip
Our notation and terminology is intended to be standard, or have a clear
meaning, e.g., $o.t.$ for order type, $card$ for cardinality.  A catalogue of
possible exceptions follows.  When forcing, $p \leq q$ means $q$ gives more
information.  Closed unbounded sets are \it clubs \rm.
The set of limit points of a
set $X$ of ordinals is denoted by $X^\prime$.  $A \Delta B$ is the symmetric
difference of $A$ and $B$, and $A\setminus B$ is the relative complement of
$B$ in $A$.  For ordinals, $\alpha \leq \beta$, $[\alpha,\ \beta)$  is the
half-open interval $\{\gamma: \alpha \leq \gamma < \beta\}$.
The notation for the three other
intervals are clear.  It should be clear from context whether the open
interval or the ordered pair is meant.
$OR$ is the class of all ordinals.
For infinite cardinals, $\kappa,\ H_\kappa$ is
the set of all sets hereditarily of cardinality $< \kappa$, i.e. those sets $x$
such that if $t$ is the transitive closure of $x$, then $card\ t <
\kappa$.  For ordinals $\alpha,
\ \beta$, we write $\alpha >> \beta$ to mean that $\alpha$ is MUCH greater than
$\beta$; the precise sense of how much greater we must take it to be
is supposed to be clear from context.  For models, $\Cal M,\ Sk_{\Cal M}$
denotes the Skolem operation in $\Cal M$,
where the Skolem functions are obtained
in some reasonable fixed fashion.  In this paper, we often suppress mention of
the membership relation as a relation of a model, but it is always intended
that it be one.  Thus, $(M,\ A,\ \cdots)$ denotes the same model as
$(M,\ \in,\ A,\ \cdots)$.
All other notation is introduced as needed (we hope).

\bigskip
\bigskip

\define\nar{\narrower\smallskip\noindent}
\bigskip

\subheading{\S 1.  THE FORCINGS}
\bigskip
\demo{(1.1)  Definition}  If $g$ is a function, $\overline g = \{x \in
dom\ g: g(x) = 1\}$.
\enddemo
\bigskip
\demo{(1.2) Definition}  If $\lambda$ is a infinite cardinal, $T \subseteq
\lambda$, then $g \in  S(T,\ \lambda)$ iff there's $\delta = \delta(g) \in
(\lambda,\ \lambda^+)$ such that
$g:(\lambda,\ \delta) \rightarrow \{0,\ 1\}$ and for
all $\alpha \in (\lambda,\ \delta]:$
\enddemo
{\nar
$(\ast)_{\alpha,\ g}\ (card\ \alpha)^{L[T,\ g|\alpha]} = \lambda$ (we say:  $g$
promptly collapses $\alpha$).
\medskip}

$\bold {S}(T,\ \lambda) = (S(T,\ \lambda),\ \subseteq)$.
\bigskip

\demo{(1.3) Definition} Let $\kappa$ be an infinite cardinal, $T \subseteq
\kappa^+,\ g \in S(T,\ \kappa^+)$.
\define\k{\kappa}
$\overset{\rightharpoonup}\to{b}^g = (b^g_\alpha : \alpha \in
(\k^+,\ \delta(g)])$ is a sequence of
almost disjoint cofinal subsets of successor ordinals
$\beta \in (\kappa,\ \kappa^+)$ which are multiples of $3$,
such that for all
$\alpha \in (\kappa,\ \delta(g)],\ (b^g_\xi: \xi \in (\kappa^+,\ \alpha])$ is
canonically defined in $L[T,\ g|\alpha]$.
\enddemo
\bigskip

\demo{(1.4)  Definition} With $\kappa,\ T,\ g$ as in (1.5), $p =
(\ell(p),\ r(p)) \in P_{\kappa,\ T,\ g}$ iff
\roster
\item $\ell(p) \in  S^+(T \cap \kappa,\ \kappa)$,
\item if $\alpha \in (\kappa,\ \delta(\ell(p))),\ \alpha = 3\alpha^\prime + 1$,
then $\ell(p)(\alpha) = 1$ iff $\alpha^\prime \in T$. (we say: $\ell(p)$ codes
T),

\item $r(p):dom\ r(p) \rightarrow \kappa^+,
\ dom\ r(p) \in
[dom\ g]^{<\kappa^+}$, and whenever
$\alpha \in  dom\ r(p)$, $r(p)(\alpha) \leq \xi \in b^g_\alpha
\cap \delta(\ell(p)),\ \ell(p)(\xi) = g(\alpha)$,

\item if $\alpha_1,\ \alpha_2 \in dom\ r(p)$ and $g(\alpha_1) \neq
g(\alpha_2)$, then $b^g_{\alpha_1} \setminus r(p)(\alpha_1) \cap
b^g_{\alpha_2} \setminus r(p)(\alpha_2) = \emptyset$.
\endroster
For $p,\ q \in  P_{\kappa,\ T,\ g}$, $p \leq q$ iff $\ell(p) \subseteq
\ell(q)$, $r(p) \subseteq r(q)$; $\bold P_{\kappa,\ T,\ g} =
(P_{\kappa,\ T,\ g},\ \leq)$.
\enddemo
\bigskip

\demo{(1.5) Definition}  Let  $\kappa$ be an infinite cardinal.
For $n \leq
\omega$, let $\kappa_n$ be $\kappa^{+n}$.  Let  $Z \subseteq \kappa_\omega$ be
such that for all $n \leq \omega$,
$H_{\kappa_n} = L_{\kappa_n}[Z \cap \kappa_n]$.
$p \in S(\kappa,\ Z) = S(\kappa)$ iff $dom\ p = \omega$, for all $n < \omega,
\ p(n) = (\ell(p(n)),\ r(p(n))),\ \ell(p(n)) \in
S(Z \cap \kappa_n,\ \kappa_n)$ and
$p(n) \in  P_{\kappa_n,\ Z \cap \kappa_{n+1},\ \ell(p(n+1))}$.
For $p$, $q \in
S(\kappa),\ p \leq q$ iff for all $n < \omega,\ \ell(p(n)) \subseteq
\ell(q(n)),\ r(p(n)) \subseteq r(q(n))$.
$\bold {S}(\kappa) = \bold {S}(\k,\ Z) =
(S(\kappa),\ \leq)$.

If $k < \omega,\ S_k(\kappa) = S_k(\k,\ Z) = \{p|[k,\ \omega): p \in
S(\kappa)\}; \leq_k$ is the obvious projection of $\leq$ onto
$S_k(\kappa)$.  $\bold {S}_k(\kappa) = \bold {S}_k(\k,\ Z) =
(S_k(\kappa),\ \leq_k)$.
\enddemo
\bigskip
\bigskip

\subheading{\S 2.  THE RESULTS}
\bigskip

Our ultimate goal in this section will be to prove that for cardinals
$\kappa$ with $\kappa \geq \aleph_2$,
for all $k < \omega,\ \bold S_k(\kappa)$ is $(\kappa_k,\ \infty)$- distributive.
As will be clear from what follows, by this we mean that the intersection
$\kappa_k$ open dense sets is dense, and not the weaker notion involving
fewer than $\kappa_k$ open dense sets.  We denote the latter notion by
$(< \kappa_k,\ \infty\text{-distributive}$.  A useful first step will be
to establish something stronger than this latter notion.

\proclaim{(2.1)  Proposition}  For all $k < \omega,\ \bold S_k(\kappa)$
is $< \kappa_k$- complete.
\endproclaim

\demo{Proof}  Let $\theta < \kappa_k,\ (p_i: i < \theta)$ be a $\leq_k$-
increasing sequence from $S_k(\kappa)$.  For $i < \theta,\ k \leq n < \omega$,
let $\delta_i(n) = \delta(\ell(p_i(n)))$, so, for such $n$,\newline
$\ (\delta_i(n): i <
\theta)$ is non-decreasing.  Let $\delta(n) = sup\ \{\delta_i(n):i <
\theta\}$.  Let $\ell(p(n)) = \bigcup\{\ell(p_i(n)):i < \theta\}$,
$r(p(n)) =
\bigcup\{r(p_i(n)): i < \theta\}$, and let
$p(n) = (\ell(p(n)),\  r(p(n)))$,
for $k \leq n < \omega$.  We shall prove that $p \in  S_k(\kappa)$.  The
only difficulty is to prove that for
$k \leq n < \omega,
\ (card\ \delta(n))^{L[Z \cap \kappa_n,\ \ell(p(n))]} = \kappa_n$.
If $\theta$ is a successor
ordinal or $\delta(n) = \delta_i(n)$ for some $i < \theta$, this is clear.
Otherwise, $\delta(n)$ is a limit ordinal of cofinality $\leq\ cf\ \theta <
\kappa_n$, so, by the Covering Lemma, already (cf $\delta(n))^L < \kappa_n$.
But then, since \newline $(\forall \alpha <
\delta(n))(card\ \alpha)^{L[Z \cap \kappa_n,\ \ell(p(n))]}
\leq \kappa_n$, the conclusion is
clear.
\enddemo
\bigskip

\noindent(2.2)  Before proving the main lemma of the section, in (2.3),
it will be helpful to simply remark (the proofs are easy, and the reader may
consult \cite{2} for an outline) that letting $\dot G$  be the canonical name
for the generic, letting $\dot {G}(n)$ be the canonical name for
$\{\ell(p(n)): p \in \dot G\}$, and letting $\dot {W}(n)$ be the
canonical name for $\bigcup \dot {G}(n)$, then for all $k < \omega$,
\medskip

{\nar $\Vdash_{\bold {S}_k(\kappa)}$
\lq\lq $(\forall k \leq n < \omega)\dot {W}(n),\ Z \cap
\kappa_n \in L[Z \cap \kappa_k,\ \dot {W}(k)]$ ". \medskip}

\noindent
We shall use a variant of this fact with no further comment below, in the
proof of the main lemma.  We note only that by an easy density argument, it
can be shown that for $k \leq n < \omega$ and
$\alpha \in [\kappa_{n + 1},\ \kappa_{n + 2})$,
there is $\eta < \kappa_{n + 1}$ such
that whenever $\xi \in b^{\dot {W}(n+1)|\alpha}_\alpha \setminus \eta,
\ \dot {W}(n)(\xi)= 0 \Rightarrow \dot {W}(n)(\xi + 1) =
\dot {W}(n + 1)(\alpha)$, and that $\{\xi \in
b^{\dot {W}(n + 1)|\alpha}_\alpha: \dot {W}(n)(\xi) = 0\}$ is cofinal in
$\kappa_{n + 1}$.  Thus, $\dot {W}(n + 1)(\alpha)$ is read by:
$\dot {W}(n + 1)(\alpha) = i$ iff there is a final segment $x \subseteq
b^{\dot {W}(n + 1)|\alpha}_\alpha$ such that for all $\xi \in x,\
\dot {W}(n)(\xi) = i$.
\bigskip

\noindent(2.3) We are now ready for the main Lemma.

\proclaim{Lemma}  For all $k < \omega$, $\bold S_k(\kappa)$ is
$(\kappa_k,\ \infty)$-distributive.
\endproclaim

\demo{Proof}  We first note that it suffices to prove that for all $k <
\omega$

{\nar
$(\ast)_k$:
Let $p_0 \in S_k(\kappa)$,
let $\chi$ be regular $\chi >>
2^{2^{\kappa_\omega}}$; let
$<(*)$ be a well-ordering of $H_\chi$ in type $\chi$; let
$\Cal M = (H_\chi,\ <(*),\ \{\bold S_k(\kappa)\},\ \{ Z\},\ \{p_0\})$;
let $\Cal N \prec \Cal M,\ \kappa_k + 1 \subseteq |\Cal N|,
\ card\ |\Cal N| = \kappa_k$.  Then there
is $p_0 \leq_k p^*$ which is
$(\Cal N,\ \bold S_k(\kappa) \cap |\Cal N|)$-generic.
\medskip}

\noindent
The argument that $(\ast)_k$ suffices is well-known, so fix the above data.
Without loss of generality,
we may assume that $[|\Cal N|]^{<\kappa_k} \subseteq |\Cal
N|$.  It will often be convenient to work with the transitive collapse of
$\Cal N$, so let $\pi:\hat{\Cal N}\rightarrow \Cal N$ be the inverse of the
transitive collapse map; thus, $[|\hat{\Cal N}|]^{<\kappa_k} \subseteq
|\hat{\Cal N}|$.  Let $\sigma = \pi^{-1} =\ $
the transitive collapse map.
If $X \subseteq |\hat{\Cal N}|$ and $(\hat{\Cal N},\ X)$
is amenable, then we
let $\pi(X) = \bigcup\{\pi(a\cap X):a\in |\hat{\Cal N}|\}$, and
similarly for $\sigma(Y)$ if $(\Cal N,\ Y)$ is amenable.  We let
$\hat {\k}_n=\sigma(\kappa_n)$.  We also let $\theta_n = sup\ (|\Cal
N| \cap \kappa_n)$.

For $k < n < \omega$, note that $\hat {\k}_n=
(\kappa_n)^{\hat{\Cal N}}$, and that
$\hat {\k}_{k + 1}= \theta_{k + 1}$.
Note that by applying Proposition 2.1 to
forcing over $\hat{\Cal N}$ with $\sigma(\bold S_{k + 1}(\kappa))$, we easily
construct $\hat p \in  \sigma(\bold S_{k + 1}(\kappa))$ which is
$(\hat{\Cal N},
\ \sigma(\bold S_{k + 1}(\kappa)))$-generic, such that
$\sigma(p_0)|[k + 1,\ \omega)$ is extended by
$\hat p$, in $\sigma(\leq_{k + 1})$,
such that for $k + 1 \leq n < \omega$, $\hat p(n) \subseteq|\hat{\Cal N}|$
and all proper initial segments of $\hat p(n)$ lie in $|\hat{\Cal N}|$.
In view of the discussion in (2.2), for forcing over $\hat{\Cal N}$,

\medskip
{\nar
$\hat{\Cal N}[\hat p] \models$
\lq\lq $(\forall n)(k + 1 \leq n < \omega \Rightarrow
\hat p(n) \in
L[\sigma(Z \cap \kappa_{k + 1}),\ \hat p(k + 1)]$".\medskip}

\noindent
Thus, $\hat{\Cal N}[\hat p] \models$
\lq\lq $\sigma(Z \cap \kappa_\omega) \in
L[\sigma(Z \cap \kappa_{k + 1}),\ \hat p(k + 1)]$".

A crucial observation is:

\proclaim{Proposition} $OR \cap |\hat{\Cal N}| < ((\hat {\k}_{k+1})^+)^L$.
\endproclaim

\demo{Proof}  Let $\hat\theta = OR \cap |\hat{\Cal N}|,\  \theta =
sup\ (OR \cap |\Cal N|)$.  Note that \newline
$\pi|L^{\hat{\Cal N}}:L_{\hat\theta}\rightarrow_{\Sigma_1}L_\theta$,
with critical point $\hat {\k}_{k+1}$.
If $\hat\theta \geq ((\hat {\k}_{k+1})^+)^L$, then $0^\sharp$ exists,
which proves the Proposition.
\enddemo
\medskip

\noindent
Thus, $(cf\ \hat {\k}_n)^L \leq (cf\ \hat {\k}_{k+1})^L$, for all $k + 1 \leq
n < \omega$.  Typically, of course, $\hat {\k}_{k + 1}$ is a
$(\text {regular cardinal})^L$.  Let
$x_{k + 1} = Z \cap \hat {\k}_{k + 1}$,
$h_{k + 1}= \ell(\hat {p}(k + 1))$.

We shall construct
in $V,\ \hat {p}(k)$ which is $|\hat{\Cal N}|$-generic for $\bold
{P}_{\kappa_k,\ x_{k + 1},\ h_{k + 1}}$, as defined in $\hat{\Cal N}$.
Among other properties, letting $h_k = \ell(\hat p(k))$,
$h_k$ will code $h_{k + 1}$.  This will be \bf clear \rm from the
construction;
we shall use this fact \bf before \rm showing that
$(cf\ \hat {\k}_{k + 1})^{L[Z \cap \k_k,\ h_k]} = \kappa_k$.
This is exactly what is required to show
that if we define $p$ by letting $p(n) = (\pi(\ell(\hat p(n))),
\ \pi(r(\hat p(n))))$
(recall our convention about $\pi(X)$ for $(\hat{\Cal N},\ X)$
amenable), then $p
\in S_k(\kappa)$ (and $p$ is $|\Cal N|$-generic for $\bold S_k(\kappa)
\cap |\Cal N|)$.

We shall have
$h_k = \ell(q_{\kappa_k})$, $r(\hat p(k)) = r(q_{\kappa_k})$,
where $q_i = (\ell(q_i),\ r(q_i))$ and $(q_i: i \leq \kappa_k)$
is defined recursively in $V$, with $q_0 =
\sigma(p_0(k))$.  For this, in $V$, we let $(D_i: i < \kappa_k)$ enumerate
the dense subsets, in $|\hat{\Cal N}|$, of
$\bold P_{\kappa_{k},\ x_{k + 1},
\ h_{k + 1}}$, as defined in $\hat {\Cal N}$.
For all $\theta < \kappa_k$,
$(D_i: i < \theta) \in  |\hat{\Cal N}|$, in virtue of the closure
property we have assumed for $|\hat{\Cal N}|$.
For all $i <
\theta$, we'll have $q_i \in |\hat{\Cal N}|$,
so, by the same observation, for $\theta < \kappa_k$,
$(q_i: i < \theta) \in |\hat{\Cal N}|$.

Also, for $j < \hat {\k}_{k + 1}$, letting $D(j)$ be the subset of
$\bold P_{\kappa_{k},\ x_{k + 1}, \ h_{k + 1}}$ consisting of those
$r$ with $\delta(\ell(r)) \geq j$, as defined in $\hat {\Cal N}$,
clearly $D(j)$ is dense and so is
among the $D_i$.  This will
guarantee that $sup\ \{\delta(\ell(q_i)): i < \kappa_k\} = \hat {\k}_{k + 1}$,
provided that we know that $q_{i + 1} \in D_i$.  This will be part of the
construction and will also guarantee the genericity of $p$.

For $i < \kappa_k$, we'll set $\alpha_i = \delta(\ell(q_i))$.  For limit
$\theta \leq \kappa_k$, we let
$\ell(q_\theta) = \bigcup\{\ell(q_i):i < \theta\}$,
$r(q_\theta) = \bigcup\{r(q_i):i < \theta\}$.  If $\theta < \kappa_k$,
by the covering argument of
the proposition of (2.1), these are always conditions, and, if $\theta <
\kappa_k$, as noted above, $(q_i: i < \theta) \in |\hat{\Cal N}|$,
so also
$q_\theta \in |\hat{\Cal N}|$.  So, we must define $q_{i+1}$, where our
crucial work is done.

For each $\alpha_i \leq \alpha < \gamma < \hat {\k}_{k + 1}$,
$\alpha$ a limit ordinal, we define
$p^{\gamma,\ \alpha,\ 1} \geq q_i$ as follows:
$r(p^{\gamma,\ \alpha,\ 1}) = r(q_i)$; if
$\alpha_i \leq \beta < \gamma$ and $\beta \equiv 1\ (mod\ 3)$ then \newline
$\ell(p^{\gamma,\ \alpha,\ 1})(\beta) = 0$
if $\beta^\prime \not\in Z\ \&\ = 1$, if $\beta^\prime \in Z$, where
$\beta^\prime$ is such that $\beta = 3\beta^\prime + 1$.  If $\gamma \geq
\alpha_i + \kappa_k$, we fix a subset
$b \in |\hat{\Cal N}| \cap L,\ b \subseteq \k_k$
which codes a well-ordering of $\kappa_k$ in type
$\gamma$, and for $\beta < \kappa_k$, we set
$\ell(p^{\gamma,\ \alpha,\ 1})(\alpha_i + 3\beta + 2) = 0$ if
$\beta \not\in b\ \&\ =\ 1$ if $\beta \in b$.
If $\alpha_i + \kappa_k \leq \beta < \gamma$
and $\beta \equiv 2\ (mod\ 3)$, we set
$\ell(p^{\gamma,\ \alpha,\ 1})(\beta) = 0$.
Similarly, if  $\gamma < \alpha_i + \kappa_k$, we set
$\ell(p^{\gamma,\ \alpha,\ 1})(\beta) = 0$
for all $\alpha_i \leq \beta < \gamma$ such that $\beta \equiv 2\ (mod\ 3)$.

If $\alpha_i \leq \beta < \gamma$ and for some $\tau \in dom\ r(q_i),
\ \beta \in b^{h_{k+1}}_\tau \setminus r(q_i)(\tau)$, then
$\ell(p^{\gamma,\ \alpha,\ 1})(\beta) = h_{k+1}(\tau)$.  Note that in virtue
of (4) of (1.4), this is well-defined.  For all other successor ordinals,
$\alpha_i \leq \beta \gamma$ which are multiples of $3$, we set
$\ell(p^{\gamma,\ \alpha,\ 1})(\beta) = 0$.

Now, suppose $\beta$ is a limit ordinal, $\alpha_i \leq \beta < \gamma$.
We set \newline
$\ell(p^{\gamma,\ \alpha,\ 1})(\beta) = 0$,
\bf unless \rm $\beta = \alpha\ \&\ =\ 1)$, if $\beta = \alpha$ (in
this case, we mark $\alpha$).

Then, let $p^{\gamma,\ \alpha,\ 2} \geq p^{\gamma,\ \alpha,\ 1}$
be chosen canonically in $D_i$.  Now
$(\gamma,\ \alpha) \mapsto p^{\gamma,\ \alpha,\ 2}$
is definable in $\hat{\Cal N}$, and so, for each $\gamma$,
in $\hat{\Cal N}$, we can compute a bound,
$\eta(\gamma) < \hat {\k}_{k + 1}$, for
$sup\ \{dom\ \ell(p^{\gamma,\ \alpha,\ 2}): \alpha_i \leq \alpha < \gamma,
\ \alpha$ a limit ordinal $\}$, as a function of $\gamma$.
Iterating $\eta$ in $\hat{\Cal N}$ gives us a club, $E_i$, of
$\hat {\k}_{k + 1},\ E_i
\in |\hat{\Cal N}|$.  Now,
$(H_{\hat {\k}_{k + 2}})^{\hat {\Cal N}} =
L_{\hat {\k}_{k + 2}}[\sigma(Z) \cap \hat {\k}_{k + 2}]$, so all clubs of
$\hat {\k}_{k + 1}$ which lie in $|\hat {\Cal N}|$, and, in particular, $E_i$,
lie in $L[\sigma(Z) \cap \hat {\k}_{k + 2}]$.  Already in $L,
\ card\ \hat {\k}_{k + 2} = card\ \hat {\k}_{k + 1}$.  So, in
$L[\sigma(Z) \cap \hat {\k}_{k + 2}]$ there
is $\theta < (\hat {\k}_{k + 1})^+$
such that all clubs of $\hat {\k}_{k + 1}$ which
lie in $|\hat {\Cal N}|$, in fact, lie in
$L_\theta[\sigma(Z) \cap \hat {\k}_{k + 2}]$.  This, however, readily gives us
that unless $(card\ \hat {\k}_{k + 1})^{L[\sigma(Z) \cap \hat {\k}_{k + 2}]} =
\kappa_k$ (and in this case, there is no problem in proving that
$q_{\kappa_k}$ is a condition),
there is a club $C$ of $\hat {\k}_{k + 1}$, $C \in
L[\sigma(Z) \cap \hat {\k}_{k + 2}]$, such that $C$ grows faster than any club
of $\hat {\k}_{k + 1}$ which lies in $|\hat {\Cal N}|$.  In particular,
$C$ grows faster than $E_i$, so that for sufficiently large
$\gamma < \hat {\k}_{k + 1}$, all $E_i$-intervals above
$\gamma$ miss $C$.  In $V$,
fix $C^* \subseteq C,\ o.t.\ C^* = \kappa_k,\ C^*$
a club of $\hat {\k}_{k + 1}$.

The idea of the above is that in constructing $p^{\gamma,\ \alpha,\ 1}$, we
have \lq\lq marked" $\alpha$ and our hope is that in passing from
$p^{\gamma,\ \alpha,\ 1}$ to $p^{\gamma,\ \alpha,\ 2}$, we have not
inadvertently \lq\lq marked" anything else.  While this is too much to hope
for, in general, we shall be able to get
that we have not marked anything else \bf in \rm $\bold C$,
provided we choose $\gamma$ sufficiently large so that every interval of
$E_i$, above $\gamma$, misses $C$.  So, GOOD's winning strategy, finally, to go
from $q_i$ to $q_{i + 1}$, is to take $\gamma$ to be the least ordinal
$> \alpha_i,\ \gamma \in C$ which, as above, is sufficiently large that
the interval $[\gamma,\ \eta(\gamma)) \cap C = \emptyset$, \bf and such that
there is \rm $\bold {\alpha^* \in [\alpha_i,\ \gamma) \cap C^*}$
and then to take $q_{i + 1} = p^{\gamma,\ \alpha^*,\ 2}$.  Thus, GOOD has
\lq\lq marked" a member of $C^*$ and nothing else in $C$, while obtaining
$q_{i + 1} \in D_i$.

Now, since, as remarked above, we know from the construction that $h_k$ codes
$h_{k + 1}$, in $L[Z \cap \kappa_k,\ h_k]$, we can recover $\sigma(Z) \cap
\hat {\k}_{k + 2}$, and therefore $C$.  But then, by the construction, we have
that $\{\alpha \in C: (h_k(\alpha),\ h_k(\alpha + 1)) = (1,\ 1)\}$ is a
cofinal subset of $C^*$.  Thus, as required,
$(cf\ \hat {\k}_{k + 1})^{L[Z \cap \kappa_k,\ h_k]} = \kappa_k$.
This completes the proof.
\enddemo
\bigskip

\noindent (2.4)  Taken together, (2.1) - (2.3) give us the
following Lemma, which, in turn, gives us the Theorem
of the Introduction:

\proclaim{Lemma} Forcing with $\bold S(\kappa)$
preserves cofinalities, $GCH$, and if $G$ is $V\text{-generic}$ for
$\bold {S}(\kappa)$, then, in $V[G]$ there is
$W \subseteq \kappa^+$
such that $V[G] = V[W],
\ Z \in L[W,\ Z \cap \k]$ and
for all $n \leq \omega,
\ H_{\kappa_n} = L_{\kappa_n}[W]$ and for
$\kappa < \alpha < \kappa^+,
\ (card\ \alpha)^{L[W \cap \alpha]} = \kappa$.
\endproclaim

\demo{Proof}  Of course $W = \bigcup\{\overline {\ell(p(0))}:
p \in \dot G\}$.
It is a routine generalization of arguments from Chapter 1 of
\cite{1} to see that for all $k$, there is $\bold Q_k \in  V^{\bold
S_k(\kappa)}$ such that $\bold {S}(\kappa) \cong
\bold S_k(\kappa) \ast \bold Q_k$,
and $\Vdash_{\bold {S}_k(\kappa)}$ \lq\lq $\bold {Q}_k$ is
$\kappa_{k+1} -$ c.c. and $card\ Q_k = \k_{k+1}$".
Further, for $k = 0$, (2.3) gives us that $\bold S(\kappa)$ is
$(\kappa,\ \infty)$-distributive and clearly $card\ S(\kappa) =
\kappa_\omega^+$.  Thus, preservation of $GCH$ is clear, as is the
preservation of all cardinals except possibly $\kappa_\omega^+$.  The
argument here is routine: if this failed, then letting $\gamma =
(cf\ \kappa_\omega^+)^{V^{\bold S(\kappa)}}$,
for some $0 < k < \omega,\ \gamma = \kappa_k$.  But then, since
$(cf\ \kappa_\omega^+)^{V^{\bold S_k(\kappa)}} > \kappa_k$,
forcing with $\bold Q_k$ over $V^{\bold S_k(\kappa)}$
would have to collapse a cardinal $\geq \kappa_{k+1}$ which is impossible.
\enddemo
\bigskip

\demo{(2.5) Remarks and Applications}
\medskip
\roster
\item  If we start from an arbitrary $Z^\prime \subseteq \kappa_\omega$, we
can, of course, code $Z^\prime$ by first coding $Z^\prime$ into a $Z$, as
above (e.g., by coding $Z^\prime$ into $Z$ on odd ordinals), and then
proceeding as above.
\medskip
\item  In work in progress, we are attempting to develop a
combinatorial approach to coding the universe by a real (when $0^\sharp$ does
not exist).  Part of our approach is to
use the Easton product of the
$\bold {S}(\kappa)$, for $\kappa = \aleph_2$, or $\kappa$ a limit cardinal,
as a preliminary forcing, to simplify the universe before doing the main
coding.
\medskip
\item  Several people have observed that the $\bold {S}(\kappa)$ afford a
method of coding any \bf set \rm of ordinals using a \bf set \rm forcing
over models of $GCH$ where $0^\sharp$ does not exist.
This can be done as follows.  Let $X \subseteq \lambda$, and assume, without
loss of generality, that $\lambda \geq \aleph_2$.  Code $X$ into a
$Z \subseteq \lambda^{+\omega}$, where $Z$ has the properties assumed above.
Then, force with $\bold {S}(\lambda)$ to get $W$, as above.  Finally, since
$W$ reshapes the interval $(\lambda,\ \lambda^+)$, we can continue to code
$W$ down to a real, using one of the usual methods of coding by a real.
\medskip
\item  Woodin has introduced the following abstract notion of condensation.
$A \subseteq \delta$ \bf has condensation \rm iff there's an algebra,
$\Cal A \in V$ with underlying set $\delta$, such that for any generic
extension $V^\prime$ of $V$:
\medskip

{\nar (*) if $X \subseteq \delta$ and $X$ is the underlying set of a subalgebra
of $\Cal A$, and $\pi:\ ({\Cal A}^*,\ A^*) \rightarrow (\Cal A|X,\ A \cap X)$,
where $\pi$ is the inverse of the transitive collapse map,
then $A^* \in V$. \medskip }

\noindent
$\delta$ has condensation iff for all $A \subseteq \delta,\ A$ has
condensation.  This notion has been investigated by Woodin's student, D. Law,
in his dissertation \cite{3}, and by Woodin himself.

S. Friedman has observed that using (3), above, it can be shown that for any
cardinal $\mu$, we can force condensation for $\mu$ via a set forcing.
We omit the proof, except to say that according to Friedman, this is not
a routine consequence of the usual sort of condensation for $L[r]$, but rather
involves a closer look at the coding apparatus.
\endroster
\enddemo

\enddocument